\documentclass[a4paper,12pt]{amsart}
\usepackage{amssymb}
\usepackage{float}
\restylefloat{figure,table}
\usepackage[colorlinks=true]{hyperref}
\hypersetup{allcolors=[rgb]{0.9,0.1,0.4}}
\usepackage{mathrsfs}
\usepackage[sort&compress,numbers]{natbib}
\usepackage{tikz}
\usetikzlibrary{decorations.markings,%
  matrix,%
  snakes}
\usepackage{wrapfig}

\newcommand{\B}{\mathcal{B}}
\newcommand{\diagram}[3]{\matrix (#1) [matrix of math nodes,row
  sep={#2},column sep={#3},text height=1.5ex,text depth=0.25ex]}
\renewcommand{\H}{\mathscr{H}}
\newcommand{\hyp}{\H\mathrm{yp}}

\renewcommand{\P}{\mathscr{P}}
\renewcommand{\S}{\mathscr{S}}
\newcommand{\sets}{\mathrm{Sets}}
\newcommand{\X}{\mathcal{X}}

\newcommand\irregularcircle[2]{% radius, irregularity
  \pgfextra {\pgfmathsetmacro\len{(#1)+rand*(#2)}}
  +(0:\len pt)
  \foreach \a in {10,20,...,350}{
    \pgfextra {\pgfmathsetmacro\len{(#1)+rand*(#2)}}
    -- +(\a:\len pt)
  } -- cycle
}

\theoremstyle{definition}
\newtheorem*{definition}{Definition}
\newtheorem*{note}{Note}

\title{On the Philosophy of Higher Structures}
\author{Nils A.\ Baas}
\address{Department of Mathematical Sciences, NTNU, N-7491 Trondheim,
  Norway}
\email{nils.baas@ntnu.no}
\date{March 11, 2018}

\begin{document}

\maketitle

\section{Introduction}
\label{sec:intro}

Higher structures occur and play an important role in all sciences and
their applications.  In a series of papers \cite{B1, B2, B3, B4, B5,
  B6, B7, B8, B9, B10, B11, B12, B13, MHS} we have developed a framework
called Hyperstructures for describing and working with higher
structures.  The purpose of this paper is to describe and elaborate
the philosophical ideas behind hyperstructures and structure formation
in general and emphasize the key ideas of the Hyperstructure Program.

In \citet{B10} we formulated a very general principal for forming
higher structures that we now call Hyperstructures and abbreviate to
$\H$-structures. There are six basic stages in this principle of
forming $\H$-structures.

\section{The $\H$-Principle}
\label{sec:Hprin}

\begin{itemize}
\item[(I)] \emph{Observation and Detection.}\\

  \noindent Given a collection of objects that we want to study and
  give a structure. First we observe the objects and detect their
  properties, states, etc.  This is the \emph{semantic} part of the
  process. Finally we may also select special objects.\\

\item[(II)] \emph{Binding.}\\

  \noindent A procedure to produce new objects from collections of old
  objects by ``binding'' them in some way. This is the
  \emph{syntactic} part of the process.\\
 
\item[(III)] \emph{Levels.}\\

  \noindent Iterating the described process in the following way:
  forming bonds of bonds and --- important! --- using the detected and
  observed properties at one level in forming the next level. This is
  iteration in a new context and not a recursive procedure. It
  combines syntax and semantics in forming a new level. Connections
  between levels are given by specifying how to dissolve a bond into
  lower level objects. When bonds have been formed to constitute a new
  level, observation and detection is like finding ``emergent
  properties'' of the process.\\
\end{itemize}

These three steps are the most important ones, but we include three
more in the general principle.

\begin{itemize}
\item[(IV)] \emph{Local to global.}\\

  \noindent Describing a procedure of how to move from the bottom
  (local) level through the intermediate levels to the top (global)
  level with respect to general properties and states. The importance
  of the level structure lies in the possibility of manipulating the
  systems levelwise in order to achieve a desired global goal or
  state. This can be done using ``globalizers'' --- an extension of
  sections in sheaves on Grothendieck sites (see \citet{B2}).\\

\item[(V)] \emph{Composition.}\\

  \noindent A way to produce new bonds from old ones. This means that
  we can compose and produce new bonds on a given level, by ``gluing''
  (suitably interpreted) at lower levels. The rules may vary and be
  flexible due to the relevant context.\\

\item[(VI)] \emph{Installation.}\\

  \noindent Putting a level structure satisfying I--V on a set or
  collection of objects in order to perform an analysis, synthesis or
  construction in order to achieve a given goal. The objects to be
  studied may be introduced as bonds (top or bottom) in a level
  structure.\\

  \noindent\emph{Synthesis:} The given collection is embedded in the
  bottom level.\\

  \noindent\emph{Analysis:} The given collection is embedded in the
  top level.\\
\end{itemize}

Synthesis facilitates local to global processes and dually, analysis
facilitates global to local processes by defining localizers dual to
globalizers, see \cite{B1}.

The steps I--VI are the basic ingredients of what we call the
\emph{Hyperstructure Principle} or in short the
\emph{$\H$-principle}. (Corresponding to ``The General Principle'' in
\citet{B10}.)  In our opinion it reflects the basic way in which we
make or construct things. This applies to mathematics, engineering,
societies and organizations. In many ways it reflects ``The Structure
of Everything!''\\

As a generic example we may think of a collection of individuals.
\begin{itemize}
\item[(a)] We observe them and detect their qualities, properties or
  states.\\
\item[(b)] We use their ``properties'' to put them into groups
  interacting to achieve specific goals.\\
\item[(c)] We observe the groups, and bind them to groups of groups
  for specific reasons.\\
\item[(d)] We introduce mechanisms (elections, communcations,\ldots)
  making it possible to pass from local states and properties to
  global ones through the levels. For example in democracies.\\
\item[(e)] Given a collection of individuals, installation is the
  process of organizing them in an $\H$-structure including (f).\\
\item[(f)] Composition means that at any level overlapping groups may
  be turned into new groups.\\ 
\end{itemize}

Evolution is a fundamental $\H$-structured process. The way we think
also follows the $\H$-principle which really represents ``food for
thought.'' We think that extracting and formalizing the essential
parts of these important structures and processes is very useful. How
do we formalize the $\H$-principle?

\section{Hyperstructures}
\label{sec:hyp}

We will here give the idea of how to formalize the $\H$-principle into
a mathematical framework. Details will be given in a separate
paper. Given a collection of objects and consider them as elements of
a set $X$. Let $\P(X)$ be the collection of all subsets of $X$. We may
also consider sets with structures like spaces, groups, orderings,
etc.

First, we \emph{observe} the subsets in terms of an assignment
\begin{equation*}
  \Omega \colon \P(X) \to \sets
\end{equation*}
(mathematically like a functor). This observation mechanism detects
properties and (or) states $\Omega(S)$ of a subset $S \subseteq
X$. This is \emph{detection}.

Then we consider subsets with properties, i.e.\ pairs $(S,\omega),
\omega \in \Omega(S)$. To each pair we want to assign another set
$B(S,\omega)$ --- the set of \emph{bonds} --- meaning ``mechanisms''
that can bind $S$ into some kind of unity. This means that we consider
the collection of all these:
\begin{equation*}
  \Gamma = \{ (S,\omega) \mid S \subseteq X, \omega \in \Omega(S)\}
\end{equation*}
We will elsewhere discuss how these assignments for $S$ and $S'$ with
$S' \subseteq S$, $S' \cap S = \emptyset$ or $S' \cap S \neq
\emptyset$ are related.

Form bonds of these pairs in terms of an assignment:
\begin{equation*}
  B \colon \Gamma \to \sets
\end{equation*}
--- in mathematical terms also possibly like a functor.

In step III we want to create a new level, proceeding as follows:
Giving everything considered so far an index $0$ we form a new set or
collection of objects:
\begin{equation*}
  X_1 = \{b_0 \mid b_0 \in B(S_0,\omega_0), S_0 \in \P(X_0),
  \omega_0\in \Omega_0(S_0)\}
\end{equation*}
Based on $X_1$ we proceed with a new $\Omega_1, \Gamma_1$ and $B_1$ to
form $X_2$ in the same way. $B_1$ will then represent bonds of
bonds. In this way we continue the process, but notice there are new
choices of assignments at each level. Hence we end up with a
hyperstructure of order $N-\H$, specified by:
\begin{align*}
  \X &= \{X_0,\ldots,X_N\}\\
  \Omega &= \{ \Omega_0,\ldots,\Omega_N\}\\
  \B &= \{B_0,\ldots,B_N\}\\
  \partial &= \{\partial_0,\ldots,\partial_N\}
\end{align*}
where $\partial$ consists of level connections
\begin{equation*}
  \partial_i \colon X_{i + 1} \to \P(X_i),
\end{equation*}
such that $\partial_i b_i = S_i$.

\begin{definition}
  The system $\H = (\X,\Omega,\B,\partial)$ is called a hyperstructure
  of order $N$.
\end{definition}

This is a semi-formal definition. In order to make it into a formal
mathematical definition more technical conditions shall be
added. This will be done in a separate paper, \cite{MHS}.

Hyperstructures extend the ideas involved in the mathematical theory
of higher categories, see \cite{L1,L2}.

\section{Discussion}
\label{sec:disc}

Once a hyperstructure has been constructed the essence lies in the
bonds: $B_0,B_1,\ldots,B_n$. To each level of bonds we may assign new
properties or states (dynamic or static). In order to ``glue'' or put
together properties or states in a compatible way we use a
Grothendieck topology and a generalized site on $\H$ (see \citet{B2})
denoted by $J$. Then we consider the pair:
\begin{equation*}
  (\H,J) \qquad \text{($\H$ of level $n$)}.
\end{equation*}
We consider a new type of property (or state) assignment
\begin{equation*}
  (\H,J) \xrightarrow{\Lambda} \S
\end{equation*}
where $\S$ is a kind of $n$-level structure for example another
hyperstructure where possibly the levels themselves can be
hyperstructures:
\begin{equation*}
  \begin{array}{r@{\; }c@{\; }l}
    B_n &\xrightarrow{\Lambda_n} &\S_0\\
    B_{n - 1} & \xrightarrow{\phantom{\Lambda_n}} &\S_1\\
    \vdots\; && \;\vdots \\
    B_0 & \xrightarrow{\Lambda_0} &\S_n.
  \end{array}
\end{equation*}
As discussed in \citet{B2} this gives a mechanism to pass from level
$B_0$ (Bottom level) to $B_n$ (Top level) of properties or states ---
the local to global process IV in the $\H$-principle. For example
there may be a preferred state $s_0$ in $\S_0$ that we want the system
to achieve by suitable dynamical actions at the lower levels. Complex
state structures (hyperstructures) may be useful in manipulating a
system in giving a wide variety of possibilities.

For this we need a ``globalizer'' which is a sequence of assignments
\begin{equation*}
  \S_n \rightsquigarrow \S_{n - 1} \rightsquigarrow \cdots
  \rightsquigarrow \S_0
\end{equation*}
compatible with the $\partial_i$'s and the site structure $J$ (see
\citet{B2}). Often it may be easier to perturb (take dynamic actions)
at lower levels. We may here perturb at the lowest level --- states in
$\S_n$ --- and via level changes at the $\S_i$'s let the changes
propagate through the structure to achieve the desired state $s_0$.

This is a formal description of what happens often in social systems
and organizations. A formal framework may help applying this to many
more situations.

Often one may want to use a collection of objects or systems (think of
a group of individuals) to achieve a goal. In order to do so one may
have to structure (organize) the collection $X$ by putting a
hyperstructure on it --- $\H(X)$ where $X$ is embedded in e.g.\ the
top or bottom level. Then one can use a globalizer or its dual to act
on the collection $X$ to achieve the desired goal.  

This is \emph{Installation} --- we will install a hyperstructure on
$X$.  Once we have formed a collection of various hyperstructures we
may apply the $\H$-principle to form $\H$-structures of such
collections again and this goes on to any order. Our point is to show
that there is a lot of room in these new higher order universes for
forming new and interesting objects both in the abstract and physical
sense.

Installation may apply to physical, biological and abstract
systems. It may also be useful in Quantum Systems making local effects
global, and controlling global states by local manipulations.

Evolutionary processes are examples of $\H$-structure. Nature uses
object properties in forming new objects, whose properties again are
used to form the next level of objects.

Let us illustrate this by an example.

An $\H$-structure as defined in Section \ref{sec:hyp} is given by:

\begin{description}
\item[Level 0] \mbox{}
  \begin{itemize}
  \item Basic objects ($X_0$): cells
  \item Properties ($\Omega_0$): receptors (selected properties)
  \item Bonds ($B_0$): aggregates of cells formed by the selected
    receptors.\\
  \end{itemize}
\item[Level 1] \mbox{}
  \begin{itemize}
  \item Objects ($X_1$): Organs formed by aggregates.
  \item Properties ($\Omega_1$): Products made by organs (selected).
  \item Bonds ($B_1$): Aggregates of organs using the product
    properties as bonds --- forming new units.\\
  \end{itemize}
\item[Level 2] \mbox{}
  \begin{itemize}
  \item Objects ($X_2$): Individuals formed by aggregates of organs.
  \item Properties ($\Omega_2$): Various types of skills for
    selection.
  \item Bonds ($B_2$): Combination of skills of individuals forming a
    unity --- a population.
  \end{itemize}
\end{description}

This structure may be refined and extended both to higher and lower
levels. The objects of each level may be subject to selection, See
\citet{Buss1987} for a thorough discussion of how selection of units
leads to the formation of higher order (hierarchical) organizations.
Hyperstructures capture this in a general sense showing the potential
of a plethora of applications where mathematical structures will be needed.

Another illustrating example is as follows. Given a finite set $X$ of
agents and the goal is to maximize interactions in $X$. Ideally
everybody would interact with everybody realizing the complete graph
on $X$, but normally there will be constraints giving a subinteraction
graph.\\

What next?\\

Then subsets of $X$ may interact and we lift the interactions to the
power set level $\P(X)$. Lots of new interactions may occur and when
all possibilites are exhausted within the given context, we may
proceed to the next level of subsets in $\P^2(X)$ using the newly
created properties to form new bonds (subsets), etc. The potential
is enormous, and $\H$-structures lead to a lot of structural novelty.

\begin{note}
  From now on we will let the assignments of properties, phases and
  states just be called states with the understanding of this broad
  interpretation. The new thing here is that we have levels of
  observables, states, properties, etc.\ --- not just local and
  global.
\end{note}

\newpage

\section{Elaborations}
\label{sec:elab}

\subsection{$\H$-formation}
\label{sub:Hform}

Once we have constructed a hyperstructure --- basically satisfying
I--III --- it may be used as a target for state (and property)
assignments in the formation and construction of new $\H$-structures,
etc. This is a very important idea. It is similar to and extends
categorification and enrichment in forming higher categories in
mathematics.

First we use the general $\H$-principle with state and bond assignment
in $\sets$ or some other known mathematical structure, e.g.\
categories like described in \citet{B10}.  We call these first level
basic $\H$-structures and denote them by
\begin{equation*}
  \hyp_0 = \{\H(0)\}.
\end{equation*}
In this case the $\Omega_i^0$'s and $B_i^0$'s take values in
structures $\S_1^0,\ldots,\S_n^0$, and we put $\S(0) = \{\S_i^0\}$.

Then we proceed with a new set of objects --- possibly from $\hyp_0$
--- and then let the $\Omega$'s and $B$'s take values in
$\H$-structures of type $\hyp_0$: $\S_1^1,\ldots,\S_n^1$ such that
even $\S(1) = \{\S_1^1,\ldots,\S_n^1\}$ may be an
$\H(0)$-structure. The results are $\H(1)$ structures.

This gives us $\hyp_1 = \{\H(1)\}$ and in this way we proceed to form
\begin{equation*}
  \begin{array}{r@{\; }c@{\; }l}
    \hyp_2 & = & \{\H(2)\}\\
           & \vdots &\\
    \hyp_n & = & \{\H(n)\}.
  \end{array}
\end{equation*}
Similarly if we are given an $\H$-structure $\H$ of some type, say,
$\hyp_k$, we may want to make new state (property) assignments and
study local to global relations.

As already described and discussed in \citet{B2}, we form a
generalized site $(\H,J)$ and give assignments
\begin{equation*}
  \begin{array}{r@{\; }c@{\; }l}
    B_n & \rightsquigarrow & \S_0\\
    B_{n - 1} & \rightsquigarrow & \S_1\\
    \vdots\;\,\\
    B_0 & \rightsquigarrow & \S_n
  \end{array}
\end{equation*}
where the $\S_i$'s are $\H$-structures at some level, even such that
$\S = \{\S_i\}$ is an $\H$-structure itself, for example of type
$\hyp_k$. In this context one may then consider globalizers (as in
\citet{B2}) relating local and global states. A globalizer will permit
levelwise dynamics or actions in order to get to a desired global
state from given local ones, as illustrated in Figure \ref{fig:glob}.

\begin{figure}[H]
  \centering
  \begin{tikzpicture}[thick]
    % Horizontal lines
    \foreach \y in {0,1.5,3,5.25,6.75}{
      \draw (0,\y) -- (8,\y);
    }

    % Triangles
    \draw (5,5.25) -- (6,6.75) -- (7,5.25);
    \draw (1,3) -- (2,4.5) -- (3,3);
    \draw (4,1.5) -- (5,3) -- (6,1.5);
    \draw (2,0) -- (3,1.5) -- (4,0);

    % Nodes - tip of triangles
    \node (T1) at (6,6.75){$\bullet$};
    \node (T2) at (2,4.5){$\bullet$};
    \node (T3) at (5,3){$\bullet$};
    \node (T4) at (3,1.5){$\bullet$};

    % Nodes - connecting nodes on horizontal lines
    \node (L1) at (3,5.25){$\bullet$};
    \node (L2) at (6,5.25){$\bullet$};
    \node (L3) at (2,3){$\bullet$};
    \node (L4) at (5,1.5){$\bullet$};
    \node (L5) at (3,0){$\bullet$};

    % Additional nodes
    \node (A1) at (2.25,6.75){$\bullet$};
    \node (A2) at (5.5,0){$\bullet$};

    % Connecting arcs - arrow tip midway
    \begin{scope}[decoration={
        markings,%
        mark=at position 0.5 with {\arrow{>}}}] 
      \draw[postaction={decorate},->] (L1) to[out=70,in=110] (L2);
      \draw[postaction={decorate},->] (T4) to[out=70,in=110] (L4);
      \draw[postaction={decorate},->] (T3) to[out=110,in=70] (L3);
      \draw[postaction={decorate},->] (A1) to[out=70,in=110] (T1);
      \draw[postaction={decorate},->] (A2) to[out=110,in=70] (L5);
    \end{scope}

    % Labels left
    \foreach \y/\ylabel in {0/0,1.5/1,3/2,5.25/{{n - 1}},6.75/n}{
      \node at (-1,\y){$B_\ylabel$};
    }
    \node at (-1,7.5){$\H:$};

    % Labels right
    \foreach \y/\ylabel in {0/n,1.5/{{n - 1}},3/{{n - 2}},5.25/1,6.75/0}{
      \node (S\ylabel) at (9,\y){$\S_\ylabel$};
    }
    \node (Sdots) at (9,4.15){$\vdots$};

    \draw[snake = coil,%
      segment aspect = 0,%
      segment amplitude = 0.5pt,%
      ->] (Sn.north) -- (S{n - 1}.south);
    \draw[snake = coil,%
      segment aspect = 0,%
      segment amplitude = 0.5pt,%
      ->] (S{n - 1}.north) -- (S{n - 2}.south);
    \draw[snake = coil,%
      segment aspect = 0,%
      segment amplitude = 0.5pt,%
      ->] (S{n - 2}.north) -- (Sdots.south);
    \draw[snake = coil,%
      segment aspect = 0,%
      segment amplitude = 0.5pt,%
      ->] ([yshift=-0.2cm] Sdots.north) -- (S1.south);
    \draw[snake = coil,%
      segment aspect = 0,%
      segment amplitude = 0.5pt,%
      ->] (S1.north) -- (S0.south);

    % Labels middle
    \begin{scope}[font=\scriptsize]
      \node[below] at (3,0){$\{s_n\}$};
      \node[above left] at (3,1.5){$\{\hat{s}_{n - 1}\}$};
      \node[below] at (5,1.5){$\{s_{n - 1}\}$};
      \node[below] at (2,3){$\{s_{n - 2}\}$};
      \node[above right] at (5,3){$\{\hat{s}_{n - 2}\}$};
      \node[above left] at (3,5.25){$\{\hat{s}_1\}$};
      \node[below] at (6,5.25){$\{s_1\}$};
      \node[above right] at (6,6.75){$\{s_0\}$};
      \node[above left] at (2.25,6.75){$\{\hat{s}_0\}$};
      \node[below] at (5.5,0){$\{\hat{s}_n\}$};
    \end{scope}

    \node at (4,4.75){$\vdots$};
  \end{tikzpicture}
  \caption{}
  \label{fig:glob}
\end{figure}

\subsection{$\H$-processes}
\label{sub:Hproc}

In many types of systems: atomic and molecular, biological,
sociological and organizational one wants to change global states by
performing local actions. $\H$-structures are useful in organizing
such actions.

Suppose we have a system $X$ of objects and an associated system or
structure of states $\S(X)$. Often a general situation is that one
wants to change the state of the system by finding a suitable action:
Given $s_0, \hat{s}_0\in \S(X)$ and construct an action
$A \colon \S(X) \to \S(X)$ changing the state. A possible way to do
this is by installing $\H$-structures: $\H(X)$ and $\H(\S(X))$ and
construct a hyperstructured action $A(\H)$ such that
$A(\H)(\hat{s}_0) = s_0$ where $A(\H) = \{A_k(\H)\}$ and in the
notation of Figure \ref{fig:glob}:
\begin{equation*}
  A_k(\H)(\hat{s}_{n - k}) = s_{n - k}.
\end{equation*}
Such a change of state could be a change of production, changing a
material property, fusion or splitting of systems, change of political
attitudes, etc.  The point is that the lower level, local actions may
be easier to perform than the global ones. A higher order state
structure allows for more varied actions and propagation through the
levels. Local to global processes in $\H$-structures are sometimes like
structured and controlled ``butterfly effects'' as seen in non-linear
chaotic systems --- amplifying small effects. See \citet{NN2009}.

\subsection{$\H$-algebras}
\label{sub:Halg}

In an $\H$-structure with bonds $\{B_0,B_1,\ldots,B_n\}$ we may define
operations or products of bonds by ``gluing.''  If $b_n$ and $b_n'$
are bonds in $B_n$ that are ``gluable'' at level $k$, then we ``glue''
them to a new bond $b_n \, \square_k^n \, b_n'$:
\begin{center}
  \begin{tikzpicture}
    \diagram{d}{2.5em}{5em}{
      b_n & b_n'\\
      b_k & b_k'\\
    };

    \path[font=\scriptsize]
      (d-1-1) edge[snake=coil,segment aspect=0,segment
      amplitude=0.5pt,->] node[midway,left]{$\partial \circ \cdots
        \circ \partial$} (d-2-1)
      (d-1-2) edge[snake=coil,segment aspect=0,segment
      amplitude=0.5pt,->] node[midway,right]{$\partial \circ \cdots
        \circ \partial$} (d-2-2)
      (d-2-1) edge[out=-30,in=-150,<->]
      node[midway,below,xshift=6.75em]{``gluable''
        (having similar parts to be identified).} (d-2-2);
  \end{tikzpicture}
\end{center}
$(\H,\{\square_k^n\})$ gives new forms of higher algebraic
structures. We have \emph{level operations} $\{\square_k^k\}$ and
\emph{interlevel operations} $\{\square_k^n\}$.

For geometric objects $X$ and $Y$ one may define a ``fusion'' product
\begin{equation*}
  X \, \square_\H \, Y
\end{equation*}
by using installed $\H$-structures on $\H(X), \H(Y)$ and
$\H(X \sqcup Y)$, see \citet{B2}.

If in an $\H$-structure we are given a bond $b_k$ binding
$\{b_{k - 1}^i\}$ the state assignments will give a levelwise
assignment via a globalizer
\begin{equation*}
  \Lambda_{k - 1}(\{b_{k - 1}^i\}) \rightsquigarrow
  \Lambda_k(b_k).
\end{equation*}
The globalizers act as generalized pairings connecting levels. In some
cases factorization algebras connect local to global observables. The
global observables may be obtained from the local ones up to
isomorphism in perturbative field theories, see
\cite{CG,Ginot,AFR}, but not in general.

Often a tensor product $\otimes_k$ may be provided in $\S_k$ we may
have
\begin{equation*}
  \Lambda_{k - 1}(\{b_{k - 1}^i\}) = \bigotimes_i
  \Lambda_{k - 1}(b_{k - 1}^i)
\end{equation*}
and sometimes when it makes sense, $\S_{k - 1} = \S_k$, like in
topological quantum field theory:
\begin{equation*}
  \Lambda_k(b_k) \in \bigotimes_i \Lambda_{k - 1} (b_{k - 1}^i).
\end{equation*}
See also \citet{B6}.

An \emph{$\H$-algebra} will be an $\H$-structure $\H$ with ``fusion''
operations $\square = \{\square_k^n\}$.  One may also add a
``globalizer'' (see \citet{B2}) and tensor-type products as just
described. The combination of a tensor product and a globalizer is a
kind of extension of a ``multilevel operad.''

\subsection{$\H$-scaffolds}
\label{sub:Hscaff}

In many situations we have systems with resources that we want to
release. This could be energy, products, human resources, etc. Often
it takes resources to get resources released. In such situations it
may be advantageous to put a suitable hyperstructure $\H(X)$ on the
system $X$. Often release of resources at the lowest level may require
small inputs, the outputs being then inputs of the next level as in
the following scheme.
\begin{table}[H]
  \centering
  \begin{tabular}{c c c@{\hspace*{2cm}} c c c}
    \multicolumn{2}{c}{$\H(X)$} && \multicolumn{3}{c}{$\S(X)$}\\
    Bonds & Resources && States & Actions & Mechanisms\\ \cline{1-2}
    \cline{4-6}
    $B_n$ & $R_n$ && $\S_0$ & $s_0 \xrightarrow{A_0} t_0$ & $M_0$\\
    $\vdots$ & $\vdots$ && $\vdots$ & $\vdots$ & $\vdots$\\
    $B_1$ & $R_1$ && $\S_{n - 1}$ & $s_{n - 1} \xrightarrow{A_{n - 1}}
                                    t_{n - 1}$ & $M_{n - 1}$\\
    $B_0$ & $R_0$ && $\S_n$ & $s_n \xrightarrow{A_n} t_n$ & $M_n$\\
  \end{tabular}
\end{table}
Then we may imagine that a small (in a precise sense) input of
resources $r_0$ will cause an action $A_n$ by the help of a mechanism
$M_n$ to change the state from $s_n$ to $t_n$, which again will cause
a larger release of resources $r_1$. This will be used as an input and
proceed through the levels upwards until a release $r_n$ is obtained.

In favourable situations with a well designed $\H$-structure we may
have:
\begin{equation*}
  r_0 \ll r_1 \ll \cdots \ll r_n.
\end{equation*}
In this way we may think of $\H(X)$ as a \emph{``scaffold''} on
$X$. The actions will take place as in Figure \ref{fig:glob}. The
mechanisms providing the actions are designed levelwise. We may think
of this in a metaphorical way: We have a dam and want to release
energy. Start with a mechanism $M_n$ drilling small holes at the
top. The released water is organized (``by bonds'') to act on another
mechanism $M_{n - 1}$ drilling bigger holes, releasing more water,
etc. In this picture the $\H$-structure acts like a scaffold.

\begin{figure}[H]
  \centering
  \begin{tikzpicture}[thick]
    \draw (0,0) rectangle (4,4);
    \node[above] at (2,4){Dam};
    \foreach \x in {0.5,1.5,2.5,3.5}{
      \draw (\x,3.5) circle(0.25cm);
    }
    \foreach \x in {1,2,3}{
      \draw (\x,2.75) circle(0.25cm);
    }
    \node at (2,2){$\vdots$};
    \foreach \x in {1.5,2.5}{
      \draw (\x,1.25) circle(0.25cm);
    }
    \draw (2,0.5) circle(0.25cm);

    \foreach \y/\ylabel in {0.5/X_n, 1.25/X_{n - 1}, 2/\vdots,
      2.75/X_1, 3.5/X_0}{
      \node at (6,\y){$\ylabel$};
    }
    \node[above] at (6,4){$\H:$};
    \draw (8,0) rectangle (10,4);

    \foreach \y/\ylabel in {0.5/\S_n, 1.25/\S_{n - 1}, 2/\vdots,
      2.75/\S_1, 3.5/\S_0}{
      \node at (9,\y){$\ylabel$};
    }

    \begin{scope}[decoration={
        markings,%
        mark=at position 0.5 with {\arrow{>}}}] 
      \draw[postaction={decorate}] (4,2) to[out=15,in=165] (8,2);
    \end{scope}
  \end{tikzpicture}
  \caption{}
  \label{fig:scaff}
\end{figure}

This is meant to illustrate how ``small'' local actions may cause
``large'' global effects on a system suitably organized into an
$\H$-structure. We see this in biological systems, societies and
organizations as well as in systems releasing various types of energy:
chemical and nuclear (fission and fusion). Hyperstructures may improve
such processes.

\subsection{$\H$-states}
\label{sub:Hstates}

States of an $\H$-structure that are levelwise connected by a
globalizer we call $\H$-states. Desired states of bonds $b_0$ at the
lowest level may be thought of as being ``protected'' by the higher
ones from external ``noise.'' On the other hand we may think of
desired states $b_n$ at the top level as being supported or created by
the lower level ones, and may also be protected by an underlying
geometric $\H$-structure of the system.

Consider $\H = \{B_0,B_1,\ldots,B_n\}$. Start with $(S_0,\omega_0)$,
$\omega_0 \in \Omega_0(S_0), b_0\in B_0(S_0,\omega_0)$. Proceed to
$b_1(\{(b_0^{j_0},\omega_0^{j_0})\},\omega_1)$, etc.
\begin{equation*}
  \begin{array}{r@{\; }c@{\; }l}
    \omega_0^{j_0} & \text{state is ``protected" by} & b_1\\
    \omega_1^{j_1} & \text{state is ``protected" by} & b_2\\
    & \vdots &\\
    \omega_{n - 1}^{j_{n - 1}} & \text{state is ``protected" by} & b_n.
  \end{array}
\end{equation*}

\begin{figure}[H]
  \centering
  \begin{tikzpicture}[thick]
    \coordinate (S) at (0,0);
    \coordinate (w01) at (-0.7,1);
    \coordinate (w02) at (-0.8,-0.65);
    \coordinate (w03) at (0.8,-0.65);
    \coordinate (w1) at (40:3.5);
    \coordinate (wn-1) at (-2.5:5.5);

    \draw[rounded corners=1mm] (S) \irregularcircle{2cm}{1mm}
    node[xshift=-1.5cm,yshift=2.25cm]{$b_1$};
    
    \begin{scope}[rounded corners=0.125mm]
      \foreach \n in {1,2,3}{
        \draw (w0\n) \irregularcircle{0.5cm}{0.125mm};
        % \node at (w0\n){$\omega_0^{j_0}$};
      }
      \node at (w01){$\omega_0^{j_0}$};
    \end{scope}

    \node at (0,0.3){$\ddots$};

    \draw[->] (w1) to[out=193,in=70] (57:2.15);
    \filldraw[fill=white] (w1) \irregularcircle{0.5cm}{0.125mm};
    \node at (w1){$\omega_1^{j_1}$};

    \node at (4,0){$\cdots$};
    \draw (wn-1) \irregularcircle{0.5cm}{0.125mm};
    \node at (wn-1){$\omega_{n - 1}^{j_{n - 1}}$};
    \draw (-20:6.5) arc(-20:20:6cm)
      node[xshift=0.75cm,yshift=-0.5cm]{$b_n$};
  \end{tikzpicture}
  \caption{}
  \label{fig:states}
\end{figure}

The local states are protected by the higher bonds which may be
realized in various ways as fields, subspaces, etc.

If we have a globalizer of states
\begin{equation*}
  \S_n \rightsquigarrow \S_{n - 1} \rightsquigarrow \cdots
  \rightsquigarrow \S_0
\end{equation*}
we may think of levels of states
\begin{equation*}
  \{s_n\} \rightsquigarrow \{s_{n - 1}\} \rightsquigarrow \cdots
  \rightsquigarrow \{s_0\}
\end{equation*}
as states of condensation, like a multilevel Bose--Einstein
condensation ($s_0$). Another interpretation is as an $\H$ organized
form of entanglement (see \citet{B8}).

Sometimes when coherence of bonds and states fail (``equations do not
hold'') higher bonds may be introduced and used in order to solve the
``frustration'' of lack of equations and introduce a new form of unity
(with proper equations). This may go on for several levels.

\section{The structure of everything}
\label{sec:every}

Suitably formulated hyperstructures cover most types of structures in
nature and science like e.g.\ trees, graphs, hypergraphs, networks,
spaces, categories, etc.\ often considered as structures with just one
level. But in nature --- in biology --- we see lots of structures
built up with several natural levels --- the effects of evolution.

The way we think and evolution have shown us the importance of
observing ``things'' and detect properties that we use in binding old
``things'' into new ``things.'' This is the key idea in forming higher
structures, and hyperstructures give us a framework in which we can
design and perform these constructions.

Our main point is that using higher structures in the sense of
hyperstructures represents an enormous unused potential for \emph{new}
kinds of design and construction. There is a lot of room in the
``space'' of structures and the $\H$-principle is a way to create new
structures by constructing $\H$-structures on sets or collection of
objects. This gives a plethora of new structures described as
hyperstructures which is useful both for construction and handling of
these new universes.

Hyperstructures represent a universal and unifying mechanism to study
both existing objects, making new objects and studying them as
well. For example, this should lead to higher structures as follows:
\begin{equation*}
  \begin{array}{l}
    \text{$\H$-states}\\
    \text{$\H$-materials}\\
    \text{$\H$-brain states}\\
    \text{$\H$-gene structure}\\
    \text{$\H$-organization, economy}\\
    \text{$\H$-language}\\
    \;\vdots
  \end{array}
\end{equation*}

In science symmetries are important. So are limits of all kinds
(``More is different,'' see \citet{B1}). We would like to add ``Higher
is different'' --- meaning that lots of new phenomena occur in higher
structures (hyperstructures) with several levels.

As pointed out in \cite{B1, B2, B3, B4, B5, B6, B7, B8, B9, B10, B11,
  B12, B13} there are many important areas of science where we think
that $\H$-structures may turn out to be very useful and important. Let
us sum up by mentioning some:\\

\begin{enumerate}
\item Evolution, Genome structure, Cancer.
\item The brain and AI systems.
\item New materials in chemistry and physics (e.g.\ molecular
  $\H$-links, high temperature superconductors).
\item Fusing systems (including nuclear fusion) and energy production
  in general.
\item Organizations, Societies, Economics, Production systems.
\item Engineering and Architecture design.
\item New higher QM-states --- extending Efimov states and GHZ-states.
  $\H$-type condensates --- extending Bose--Einstein condensates.
\item New universes of Abstract Matter for design and construction.
\item Mathematics in many ways. Extending for example higher
  categories.\\
\end{enumerate}

What more? Time will tell!  Only our imagination limits the list.

\begin{note}
  $\H$-structures are structures that I imagined and ``saw'' as a
  child, but it took a lifetime to describe what I saw!
\end{note}

\subsection*{Acknowledgements}
I would like to thank M.~Thaule for his kind technical assistance in
preparing the manuscript.

\section*{Notes on the contributor}
\begin{wrapfigure}[8]{l}{3cm}
   \centering
   \vspace*{-0.5cm}
   \includegraphics[width=3cm]{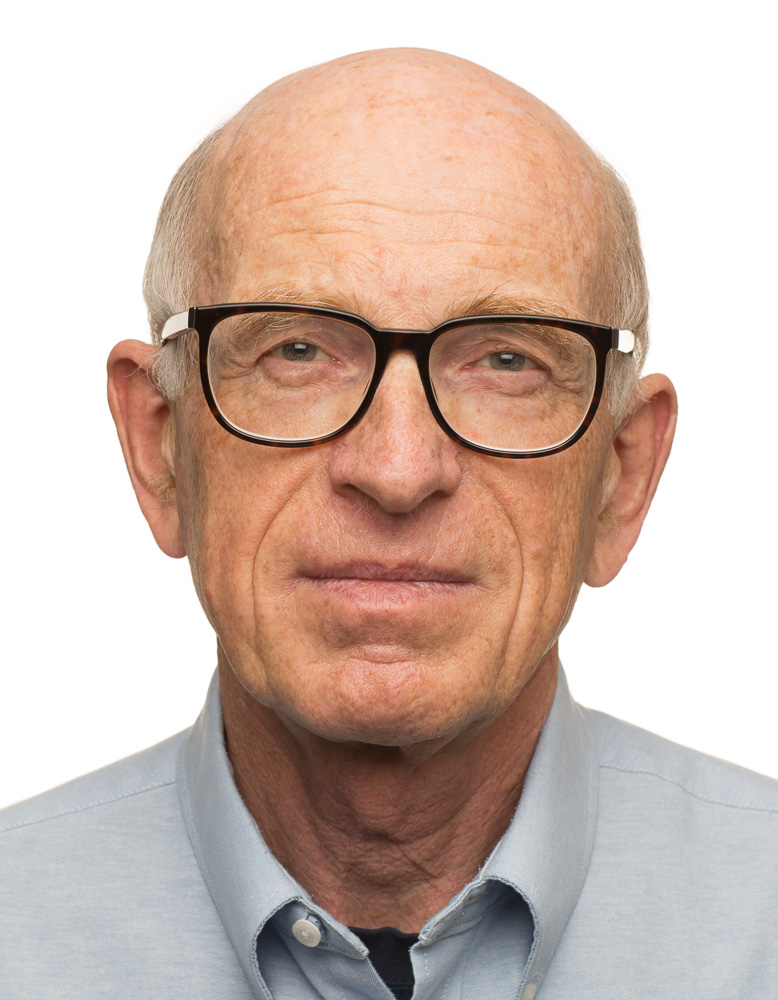}
\end{wrapfigure}

\noindent \textbf{Nils A.\ Baas} was born in Arendal, Norway, 1946.
He was educated at the University of Oslo where he got his final
degree in 1969.  Later on he studied in Aarhus and Manchester.  He was
a Visiting Assistant Professor at U. Va. Charlottesville, USA in
1971--1972.  Member of IAS, Princeton in 1972--1975 and IHES, Paris in
1975.  Associate Professor at the University of Trondheim, Norway in
1975--1977 and since 1977, Professor at the same university till date.
He conducted research visits to Berkeley in 1982--1983 and 1989--1990;
Los Alamos in 1996; Cambridge, UK in 1997, Aarhus in 2001 and 2004.
Member IAS, Princeton 2007, 2010, 2013 and 2016.  His research
interests include: algebraic topology, higher categories and
hyperstructures and topological data analysis.

%\bibliographystyle{plainnat}
%\bibliography{PHS_bib}

%
% Imported bbl
%

%
% End document
%
\end{document}